\documentclass[12pt]{amsart}
\usepackage{amsmath}
\usepackage{amssymb}

\theoremstyle{plain}
\newtheorem{lemma}{Lemma}[section]

\theoremstyle{definition}
\usepackage{hyperref}
\usepackage{float}

\usepackage{graphicx}
\usepackage{svg}
\newtheorem{thm}{Theorem}[section]

\newtheorem{lem}[thm]{Lemma}
\def\HH{\mathbb{H}}

\def\im{\mathrm{Im}\,}

\def\xx{\HH/g2}

\def\ZZ{\mathbb{Z}}
\def\CC{\mathbb{C}}
\def\RR{\mathbb{R}}
\def\QQ{\mathbb{Q}}

\def\fp{\mathbb{F}_p}

\def\gl2{\mathrm{GL}(2, \ZZ)}
\def\sl2{\mathrm{PSL}(2, \ZZ)}
\def\slr{\mathrm{PSL}(2, \RR)}
\def\g2{\Gamma(2)}

\def\xx{\HH/\g2}

\title{Eisenstein integers and equilateral ideal triangles}

 \author[McShane]{Greg McShane}
\address{Institut Fourier 100 rue des maths, BP 74, 38402 St Martin d'H\`eres cedex, France}
\email{mcshane at univ-grenoble-alpes.fr}

\thanks{This work was partially funded by the Equipe Action ToFu part of Persyval-Lab}


\begin{document}

\maketitle

\begin{abstract} 
We discuss the relationship between Penner's $\lambda$-length
and the norms of Eisenstein integers. This leads to a geometric
proof of the fact, attributed to Fermat, that every prime $p$ of the form $3k + 1$
is the norm of an Eisenstein integer that is can be written
as $a^2 - ab + b^2$ for some $a,b \in \ZZ$.
\end{abstract}

\section{Introduction}

A celebrated result of Fermat characterises
those primes which can be written as the sum of two squares

\begin{thm}[Fermat]\label{main}
Let $p$ be a prime then the equation
$$x^2 + y^2 = p $$
has a solution in integers  iff  $p =2$ or $p$ is of the form $4k + 1$.
\end{thm}

An alternative formulation of this  is that
a prime $p$ is the norm of a Gaussian integer
if and only if $p =2$ or $p$ is of the form $4k + 1$.

There are many proofs of this result see for example Elsholtz
\cite{elsholtz} for a survey of some of them.
In particular, about  1984 Heath-Brown published a proof of Theorem
\ref{main} (see also \cite{aigner2})
  in the journal of the Oxford University undergraduate mathematics society. 
His proof arose from a study of the account of Liouville’s papers on
identities for parity functions. Zagier's celebrated one line proof
\cite{zagier} is a clever reformulation of this argument but
probably the most elegant incarnation is the recent proof given by
Dolan \cite{dolan}. 
Heath-Brown studies the action of a pair of involutions on a finite
set.  
Any fixed point of one of the involutions is a solution 
to $x^2 + y^2 = p$,
the other has a single fixed point which he gives explicitly
this means that the set has an odd number of points so the first
involution must fix at least one of them.
In a previous work \cite{vlad} 
we gave a geometric proof which, in some sense,
mimics Heath-Browns method.
Our approach was based on the analysis of the action of a  Klein four group 
on complete geodesics with two ends terminating at cusps, often referred to as \textit{arcs}, of
$\lambda$-length $p$ on the surface $\xx$.
The interested reader might consult \cite{springborn1, springborn2}
for a nice account of how arcs and their $\lambda$-lengths
play a role in Diophantine approximation etc.

\subsection{Eisenstein integers}

The Eisenstein integers, $\ZZ[\omega]$,
is  the ring of integers of the cyclotomic field $\QQ(\omega)$
where $\omega$ denotes  an irrational cubic root of unity.
If $a + b\omega \in \ZZ$ is  an Eisenstein integer then its norm is
$$a^2 - ab  + b^2.$$
There is an analogue of Theorem \ref{main} in this setting:

\begin{thm}\label{eisenstein}
Let $p$ be a prime then the equation
\begin{equation}\label{eisenstein norm}
a^2 - ab + b^2 = p 
\end{equation}
has a solution in integers  iff  $p =3$ or $p$ is of the form $6k + 1$.
\end{thm}

It is easy to see why this condition is necessary since 
if we make a reduction modulo $p$ then  in $\fp$
the  equation  becomes
$$\bar{a}^2 -  \bar{a}\bar{b} +\bar{ b}^2 = 0$$
and so , for $p> 3$,  $\bar{a}/\bar{b}$ is a
cubic root of $-1$ ie an element of order 6 in the group 
$\fp^*$.
It follows that 6 divides 
the order of $\fp^*$  
which is of course equal to $p-1$.

In this note we show that this  condition is sufficient.
As before our  proof uses the action  
of a group of automorphisms of  the surface $\xx$ on
arcs of $\lambda$-length $p$.

\subsection{Farey diagram and $\lambda$-lengths}

\begin{figure}[hb]
\begin{center}
\includegraphics[scale=.5]{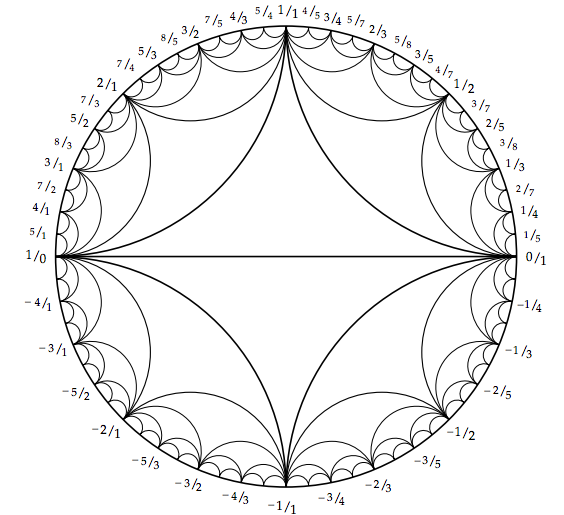} 
\end{center}
\caption{Farey diagram.}
\label{farey diagram}
\end{figure}

The  \textit{Farey diagram} (see Figure \ref{farey diagram}) is a fundamental object in the theory of Fuchsian groups. It is a tessalation of hyperbolic space by \textit{ideal triangles}.

\begin{figure}[ht]
\begin{center}
\includegraphics[scale=.4]{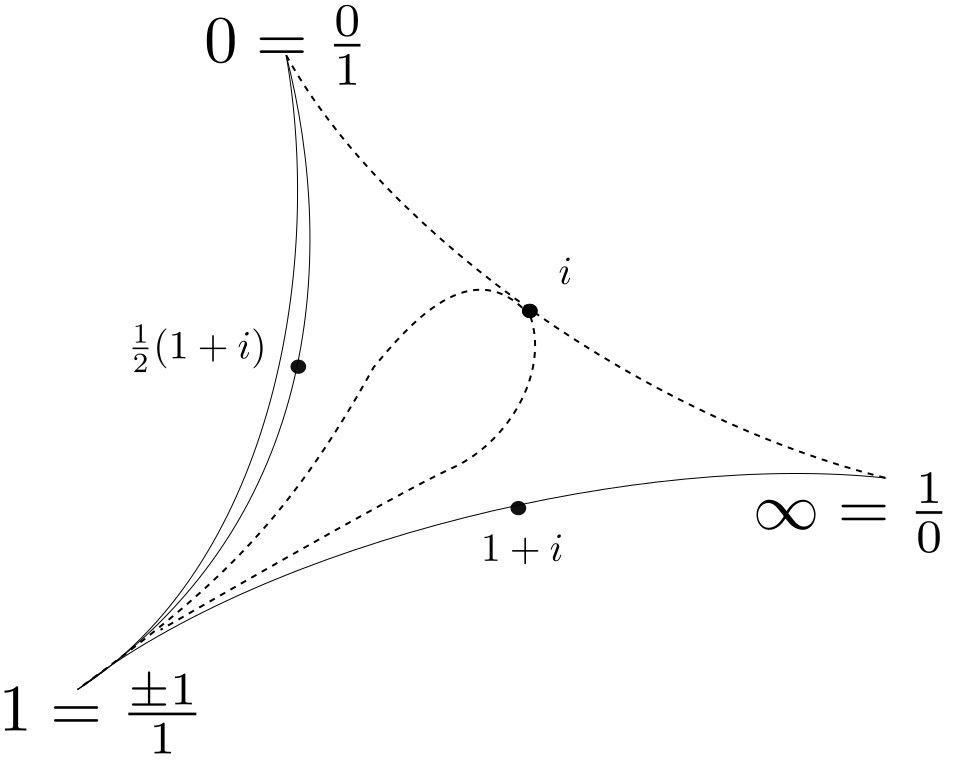} 
% \includesvg[inkscapelatex=false]{3sphere_again}
\end{center}
\caption{Three punctured sphere with cusps labelled by $\Gamma(2)$-orbits
and some arcs.
Note that the dotted geodesic in the middle is the edge 
of an ideal triangle which is embedded but not properly immersed
as two spikes meet at the cusp $\infty$.}
 \label{3punctured}
\end{figure}

The tessellation is invariant under the action of the modular group
$$\Gamma = \sl2< \slr = \mathrm{SL}(2, \RR)/\{ \pm I\}.$$
The quotient  of the upper half space by the modular group 
is the \textit{modular orbifold}. This is a non compact singular
surface with a single cusp and two cone points. 
One of the cone points lifts to the orbit of $i$  under $\Gamma$ 
and the other to the orbit of $\omega$.
The principal congruence subgroup $\g2<\Gamma$
 is a torsion free, normal subgroups of index 6
and the modular orbifold admits a degree 6 cover
corresponding to $\g2$
namely, the \textit{three punctured sphere} $\HH/\g2$ (see Figure \ref{3punctured}).
Whilst this surface is non compact it has finite hyperbolic area
equal to $2\pi$. The non compactness is due to the presence of
\textit{cusps} and in fact  $\HH/\g2$ has exactly three cusps one for each orbit of $\g2$ on the extended rationals $\QQ \cup \{ \infty \}$.

We will see (Lemma \ref{cubes}) that the reciprocal direction in 
Theorem \ref{eisenstein} can be stated in terms  of  $\lambda$-lengths
 of arcs on  this surface.
 For the purposes of this paper
 an \textit{arc}  
 is the projection to $\xx$ of a 
 Poincar\'e geodesic with both  endpoints in the 
 extended rationals $\QQ \cup \{ \infty\}$.
The geodesic edges of the Farey diagram yield a set of three such geodesic, 
 in fact the shortest arcs
 in terms of $\lambda$-length,  on the quotient surfaces.
The geodesics obtained from projecting the edges of the Farey
diagram to $\xx$
are \textit{simple}, that is they have no self intersections, 
and their complement consists of a pair of (ideal) triangles.
Of course, an arc $\alpha$ has infinite length with respect to the Poincar\'e metric  but one can define a useful geometric quantity   by truncating the geodesic ie removing a portion of infinite length  inside the cusp regions of the surface.
This idea of associating a finite length to a simple arc,
appears in Penner's work on moduli \cite{bob}
 (see paragraph \ref{lengths}.) 
on a punctured
surface to be the exponential of half the 
length of the portion outside of some fixed
system of cusp regions.
We extend this notion to non simple arcs:
Lemma \ref{calcul} shows that in the context we consider
a $\lambda$-length is always the determinant 
of an integer matrix.

\subsection{Properly immersed ideal triangles}

An \textit{ideal triangle} in $\HH$  is the convex hull 
of a triple of distinct points in $a,b,c \in \partial \HH = \RR \cup \{ \infty \}$.
We say that $a,b,c$ are the \textit{vertices of the triangle}
and that the triangle has \textit{spikes} at $a,b,c$.
We identify the  group of orientation preserving automorphisms of $\HH$
 with $\slr$ acting by linear fractional transformations.
This action is simply transitive  on
triples of distinct points in $(a,b,c) \in (\partial \HH)^3$
so there is only one ideal triangle up to isometry.

Any  ideal triangle in the Farey diagram, 
for example the one with vertices $0,-1,\infty$
projects to an embedded ideal triangle in $\xx$.
We will see later that each of the edges
of this triangle has $\lambda$-length $1$
so it is \textit{equilateral} for this notion of length.
Our triangle is also invariant by the group generated by the
hyperbolic isometry
$$ \psi: z \mapsto \frac{z+1}{-z  }.$$
The fixed points in $\CC$ of this map satisfy
$$z^2 + z + 1 = 0,$$
and there is exactly one in $\HH$ 
namely the cube root of unity $\omega= \exp(2\pi i/3)$.
The {barycenter} 
of the ideal triangle $0,-1,\infty$
is unique so it must be fixed by $\psi$ 
and so coincides with $\omega$.
Since there is a unique ideal triangle up to isometry 
the notion of barycenter is well defined for any ideal triangle and
it is the image of $\omega$ by some isometry taking $0,-1,\infty$ to
the vertices of the ideal triangle in question.

\begin{figure}[h]
\begin{center}
\includegraphics[scale=.4]{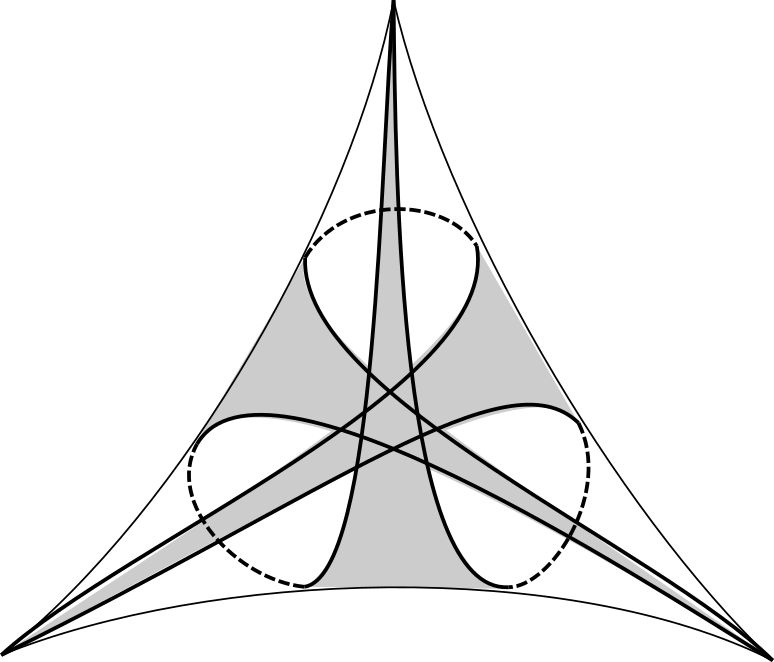} 
\end{center}
\caption{Projection of  the ideal triangle $1/3, 2/3,\infty$
to $\xx$.}
	\label{immersed triangle}
\end{figure}

Each edge of the ideal triangle $0,1/3,2/3$
has $\lambda$-length $3$ 
so again it is equilateral.
However, it does not project to  an embedded triangle
but its image is instead an immersed triangle on $\xx$
(see Figure \ref{immersed triangle}).
and further each of its spikes ends at a different cusp on 
$\xx$.
We say that such an ideal triangle,
that is one for which 
each of its spikes ends at a different cusp,
is \textit{properly immersed}.

\subsection{Sketch of proof}

We are interested in the action of a cyclic subgroup of the 
automorphisms of $\xx$ induced by $\psi$.
More precisely we consider its action on the
set of properly immersed equilateral ideal triangles in $\xx$
all of whose sides have $\lambda$-length $p$.
Let $F$ denote the cardinal of this set.
By counting incidences  (see Section \ref{counting arcs etc}) between  
arcs of $\lambda$-length $p$
and the cusps of $\xx$
we  have the following equation:
\begin{equation} \label{incidence}
3 \times F = 2 (p-2) \times 3.
\end{equation}
The RHS is $3F$ precisely because all the  triangles
are properly immersed so that each spike ends at a
different cusp on $\xx$ (see Lemma \ref{proper})
and the LHS follows from our Lemma \ref{counting}
which also  appeared in our work \cite{vlad}.
So $F = 2(p-2)$ and  if $p= 3k + 1$ then $F = 2(3k -1)$
which is not a multiple of $3$.
It follows that there are two $\psi$-orbits of length $1$
that is $\psi$ leaves two ideal triangles of this family invariant.

Then  by considering the barycenter
of a lift of one of these invariant ideal triangles
we obtain $p$ as the norm of an element of $\ZZ[\omega]$.

\subsection{The case $p=3$}

\begin{figure}[ht]
\begin{center}
\includegraphics[scale=.5]{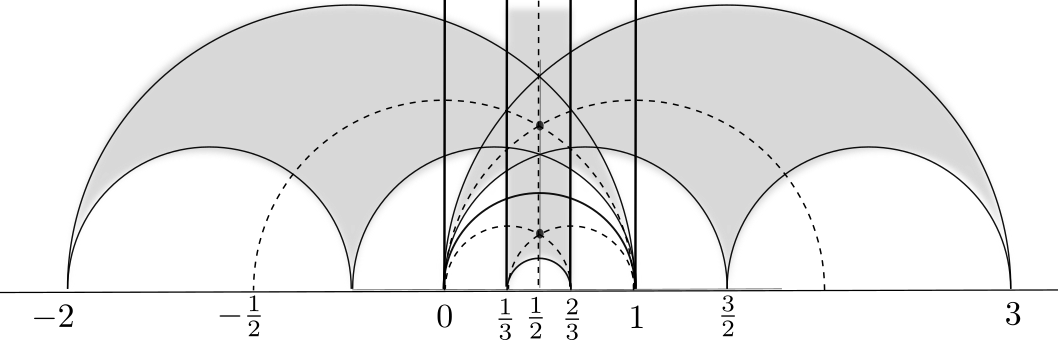} 
\end{center}
\caption{Three intersecting  lifts  of an equilateral ideal triangle.
The dotted lines intersect in two points the lower of which 
is the barycenter of $1/3, 2/3,\infty$
which coincides with the barycenter of $0, 1/2,1$.
The upper point is the barycenter of $0,1,\infty$
that is $1+\omega$.}
\end{figure}

In fact it is not difficult to show that there
are exactly two immersed equilateral 
with sides of  $\lambda$-length $3$ in $\xx$:
the  ideal triangles with vertices $0,1/3,2/3$ 
and   $0,4/3,5/3$.
Each of these is invariant under the 
automorphism induced by $\psi$
so their barycenters are in the $\sl2$ orbit 
of $\omega$.

Let us compute the barycenter of $0,1/3,2/3$ explicitly to check this. 
One has
$$3 = 2^2 - 2 + 1 = | 2\omega + 1 |^2.$$
Let
$$f :  z \mapsto \frac{z}{2z+1}$$
then the map
$$ f \circ \psi \circ f^{-1} :  z \mapsto \frac{-z+ 1}{-3z+2}$$
is conjugate to $\psi$ and permutes the vertices of our ideal
triangle as follows:
$$ \infty \mapsto 1/3 \mapsto 2/3 \mapsto \infty.$$
so it leaves our ideal triangle invariant.
One also sees from this that the barycenter of $0,1/3,2/3$
is the fixed point of $f \circ \psi \circ f^{-1}$ 
(compare with equation (\ref{eisen integers}) below)
which is just the image of $\omega$ under $f$:
$$ f(\omega) = \frac{\omega}{2\omega+1} 
= \frac{3/ 2+ \sqrt{3}i/2 }{|2 \omega + 1 |^2}
= \frac{3/ 2+ \sqrt{3}i/2 }{3}
% = \frac{3+ \sqrt{3}i }{2|2 \omega + 1 |^2}
% = \frac12 +  \frac{\sqrt{3}}{2}\frac13 i.  
$$

So we have exhibited $3$ as the norm of an Eisenstein integer namely
$2\omega + 1$. 

% In fact this is the $\lambda$-length of the arc 

\subsection{Thanks}

It is a pleasure to thank Louis Funar, Hidetoshi Masai, Robert
Penner, Vlad Sergiescu, and Xu Binbin for discussions and
suggestions over the years.

Some of the text of this paper was suggested by GitHub
Copilot\cite{copilot, vim_copilot}.

\section{Reciprocals of norms and arcs}

As we stated above an alternative formulation 
of Theorem \ref{main}  is that
a prime $p$ is the norm of a Gaussian integer
if and only if $p=2$ or $p-1$ is a multiple of $4$.
We recall the result from \cite{vlad} 
which relates $\sl2$-orbits of $i \in \CC$
to representations of a number as a sum of squares
and then discuss the analogous result for 
$\sl2$-orbits of $ \omega$.

Geometrically, we think of $p$ in Theorem \ref{main}
as the $\lambda$-length of an arc on the surface $\xx$
and the midpoint of this arc will be the projection 
of a point in the $\sl2$-orbit of $i$.
So we will begin by recalling some basic facts about
$\lambda$-lengths.

\subsection{Ford circles, $\lambda$-lengths} 
\label{lengths}

%Lemma \ref{squares} illustrates the connexion between sums of squares,
%the orbit $\sl2.i$ and Poincare geodesics. 
We now recall some standard ideas from hyperbolic geometry
necessary to define $\lambda$-length.
We define an \textit{arc} to be a Poincare geodesic
with endpoints in $\partial \HH$ a pair of extended rationals, 
that is elements of $\QQ \cup \infty$.

We denote by $F$ the set  $\{ z, \im z > 1\}$
this is a \textit{horoball in $\HH$} centered at $\infty$.
The image of $F$ under the $\sl2$ action consists of
$F$ and infinitely many disjoint discs, 
which we will refer to as \textit{Ford circles}, 
each tangent to the real line at some rational $m/n$.
We adopt the convention that $F$ is also a Ford circle of infinite radius
tangent to the extended real line at $\infty = 1/0$.
Since the interiors of the Ford circles are $\Gamma(2)$-invariant
and disjoint they project to disjoint non compact regions on $\xx$
which we will refer to as \textit{canonical cusp regions}.

\begin{figure}[ht]
\begin{center}
\includegraphics[scale=.8]{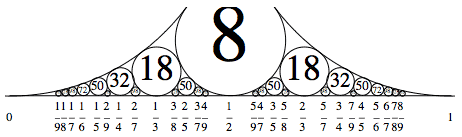} 
\end{center}
\caption{Ford circles with tangent points and curvatures.
Recall that the 
curvature of a euclidean circle is twice the reciprocal of the square of its radius.}%needs checking
\end{figure}

The following is well known \cite{ford,conway} and in any case is easily checked:

\begin{lem}\label{ford}.

The Ford circle tangent to the real line at $m/n$
has Euclidean diameter $1/n^2$.
\end{lem}

%Two Ford circles are adjacent iff their closures meet in a point
%and this point is invariably a point of the orbits $\sl2.i$

\subsubsection{$\lambda$-length}

Let $a/c, b/d$ be a pair of distinct rationals.
We define the \textit{length} of the arc 
joining these rationals 
to be half the length, with respect to the Poincaré metric on $\HH$, 
of the portion  outside of the Ford circles tangent at $a/c, b/d$ .

Following Penner \cite{bob} we define the
$\lambda$-length of an arc to be the exponential of this 
length.
 % It is a consequence  of Lemma \ref{calcul} below  that the arcs of $\lambda$-length 1 are the edges in the 
% \textit{Farey diagram} (see Figure \ref{farey diagram}).

\begin{lem}\label{calcul}
Let $a/c, b/d$ be a pair of distinct extended rationals.
Then the  $\lambda$-length of the arc joining them
 is the determinant of the matrix
$$\begin{pmatrix}
a & b \\ c & d
\end{pmatrix}.$$
Further if $a/b = 1/0$ then the arc is a vertical line 
 whose midpoint has imaginary part equal to $1/|d|$ .
\end{lem}

% It is important to note that whilst the trace of an element of
% $\slr$ is not well defined its determinant is.

\proof 
By transitivity of $\slr$ on the extended rationals we may assume 
$a/c = 1/0$ and so the determinant of the matrix is $d$.
The Ford circle at $1/0$  is $F= \{ z, \im z > 1\}$ and
at $b/d$ it is a circle of euclidean diameter $1/d^2$.
The top of this circle is tangent to the line $\{ z, \im z = 1/d^2\}$
and the $z \mapsto d^2 z $ maps this to the boundary 
of the Ford circle $\{ z, \im z \geq 1\}$.
It is easy to see that the distance between these sets,
and so the length of the segment outside the Ford circles
is $2\log d$. Thus the $\lambda$-length is equal to the 
determinant (ie $|d|$) as required.
\hfill $\Box$

Note that for any edge of the Farey diagram the determinant is $\pm
1$ so the $\lambda$-length of the arc is 1. In particular any ideal
triangle in the Farey diagram, for example the one with vertices
$0,-1,\infty$ considered above is equilateral with sides of
$\lambda$-length $1$.

\subsection{Gaussian integers}

We begin by presenting a calculation from \cite{vlad}
which relates the number of ways of writing a number as a sum of
squares that is as the norm of a Gaussian integer.

The group of integer matrices $\sl2$  acts on $\HH$ by linear fractional transformations
that is:
$$\begin{pmatrix}
a & b \\
c & d
\end{pmatrix} \in \sl2,\, z\in \HH,\, 
\begin{pmatrix}
a & b \\
c & d
\end{pmatrix}. z = \frac{az + b}{cz + d}.
$$
The key lemma that relates this  $\sl2$ action to sums of squares is:

\begin{lem} \label{squares}
Let $n$ be a positive integer.
The number of  ways of writing $n$  as a  sum of squares
$$n = c^2 + d^2$$
with $c,d$ co prime positive  integers
is equal to the number of  integers $0 < k < n-1$ co prime to $n$
such that the line
$$\{  k/n + i t,\, t >0 \}$$
contains  a point in the $\sl2$  orbit of $i$.
\end{lem}

\proof  Suppose there is such  a point which we denote  $w$.
The point $w$ is a fixed point of some  element of order 2 in $\sl2$.
Since the Ford circles are $\sl2$ invariant
this element must permute $F$ with the Ford circle tangent 
to the real line  at the real part of $w$.
So, in particular, $w$ is the midpoint of the line 
that it lies on 
and by  Lemma \ref{calcul} one has:

\begin{equation}
\label{gaussian integers}
\frac{1}{n} = \im \frac{1 }{n}(k + i)  
= \im  \frac{ai +b}{ci+d }
= \frac{\im i} {c^2 + d^2}.
\end{equation}

Conversely if $c,d$ are co prime integers 
 then there exists $a,b$ such that
 $$ad - bc = 1 \Rightarrow  
 \begin{pmatrix}
 a & b \\
 c & d
 \end{pmatrix} \in \sl2.
$$
%The real part of $w = \frac{ai +b}{ci+d }$ is
%$$ac + bd = (a,b).(c,d),$$
%and
By applying a suitable iterate of the parabolic transformation 
$z \mapsto z + 1$,
one can choose $w$ such that $0 \leq \text{Re\,} w < 1$.
So if $n = c^2 + d^2$ then $\frac{ai +b}{ci+d }$
is on one of the lines of the family in the statement.

\hfill $\Box$

\subsection{Eisenstein integers}

We replace $i$ by $\omega$ in (\ref{gaussian integers}) above:

\begin{equation}
\label{eisen  integers}
 \im  \frac{a \omega +b}{c \omega+d }
= \frac{\im \omega} {c^2  - cd  + d^2} 
= \frac{\sqrt{3}/2} {c^2  - cd  + d^2}.
\end{equation}

This expression is not so tidy but we can still exploit it 
to prove Theorem \ref{eisenstein}. 
By mimicking the proof of the previous lemma 
one can easily show:

\begin{lem} \label{cubes}
Let $n$ be a positive integer.
The number of  ways of writing $n$  as 
$$n = c^2  - cd  + d^2$$
with $c,d$ co prime positive  integers
is equal to the number of  integers $0 < k < 2n-1$ co prime to $2n$
such that the vertical line
$\{  k/2n + i t,\, t >0 \}$
contains  a point in the $\sl2$  orbit of $\omega$.
\end{lem}

The fact that the denominator is $2n$ and not $n$ 
is perhaps what is most disturbing but
 the geodesic is paired to an equilateral ideal triangle
which suits our approach:

\begin{lem} \label{cubes}
Let $p>2$ be prime and $0 < k < p$ an integer 
co prime to $2p$ such that the line $\{  k/2n + i t,\, t >0 \}$
contains  a point in the $\sl2$-orbit of $\omega$.
Then the ideal triangle 
with vertices $\infty, (k-1)/2p, (k+1)/2p$:
\begin{itemize}
\item has barycenter at $\frac{ k + i \sqrt{3}}{2n}$.
\item is equilateral  with sides of $\lambda$-length $p$.
\end{itemize}
\end{lem}

The  first part follows from the discussion above and the
second part is a corollary of  Lemma \ref{calcul} below.

\section{The three punctured sphere}

We consider $\g2$,
the principal level 2 congruence subgroup of $\sl2$.
This group acts on $\ZZ^2$, that is pairs of integers,  preserving parity.

 \begin{figure}[hb]
\begin{center}
\includegraphics[scale=.5]{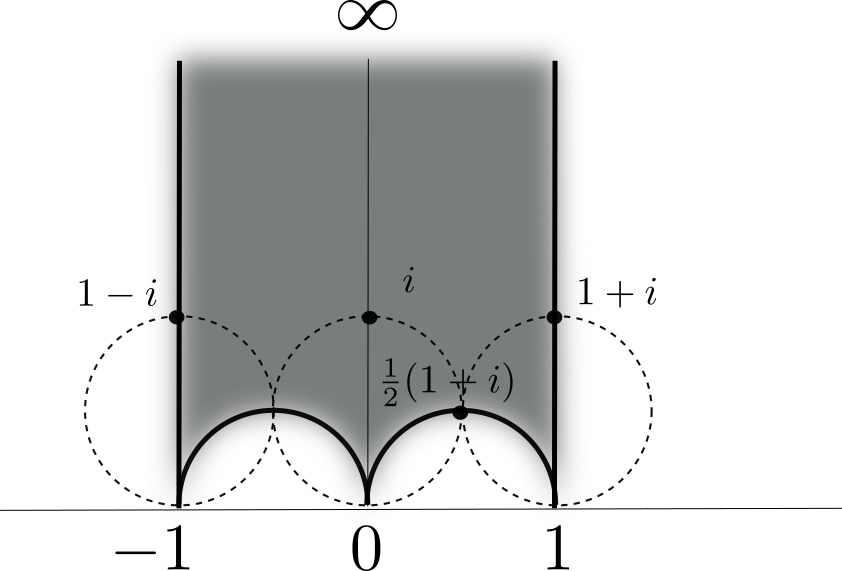} 
\end{center}
\caption{Standard fundamental domain for $\g2$ and its decomposition into ideal triangles.}
\label{fund}
\end{figure}

It also acts on $\HH$ by linear fractional transformations
and the  quotient $\xx$ is conformally equivalent to the Riemann  sphere minus three points
which we will refer to as \textit{cusps}
(see Figure \ref{3punctured}). 
Following  convention we label these cusps $0,1,\infty$ respectively
corresponding to the three $\g2$ orbits of $\QQ \cup \infty$. 
Finally, the \textit{standard fundamental domain}  for $\g2$ 
is the convex hull of the points $\infty, -1, 0 , 1$.
This region can be decomposed into two ideal triangles 
$\infty, -1, 0 $ and $ 0 , 1,\infty$ as in Figure \ref{fund}.
The edges of the ideal triangles project to three disjoint simple geodesics on $\xx$
and each edge has a \textit{midpoint} 
which is a point of the $\sl2$ orbit of $i$ (see Figure \ref{3punctured}).

\subsubsection{Cusp regions}

The image of a Ford circle on $\xx$ is a \textit{cusp region}
around one of the three cusps $0,1,\infty$.
Pairs of these cusp regions are tangent at one of the midpoints
labelled $i, 1+i, \frac12(1+i)$.
It is not difficult to see that these cusp regions 
are permuted by the automorphisms of $\xx$.
It follows that if an automorphism preserves a geodesic  joining cusps on $\xx$
then it must permute the Ford regions at each end of a lift to $\HH$.

\subsection{Automorphism groups of $\xx$}
From covering theory an isometry  of $\HH$ 
induces an automorphism of $\xx$ iff it normalises the covering group
i.e. $\g2$.
It follows that,
since $\g2$ is a normal subgroup of $\sl2$,
 the quotient group
 $$ \sl2/\g2 \simeq \mathrm{SL}(2, \mathbb{F}_2) \simeq \mathfrak{S}_3$$
acts as a group of (orientation preserving) automorphisms of the surface $\xx$.

%\subsection{Order 3 automorphisms}

In our proof of Theorem \ref{main}
it was necessary to study the action of
a Klein 4 group of automorphisms
which contained orientation reversing automorphisms.
For the proof of Theorem \ref{eisenstein}
the analysis is simpler as what is important
is  the action of a
cyclic group of order 3
generated by
$$ \begin{pmatrix}
1 & 1 \\
-1 & 0
\end{pmatrix} .$$

This automorphism:
\begin{enumerate}
\item permutes the cusps;
\item leaves each of the properly embedded  ideal triangles invariant;
\item has exactly 2 fixed points and these coincide with the
	barycenters of each of the embedded ideal triangles.
\end{enumerate}

\section{Counting arcs and ideal triangles}
\label{counting arcs etc}

Let $p> 2$ be a prime. To prove the relation (\ref{incidence})
we have to  
\begin{itemize}
\item count the number 
of arcs of $\lambda$-length $p$
incident to the cusp $1/0 = \infty$ on $\xx$.
\item 
show that all the equilateral triangles with sides of $\lambda$-length $p$ are properly immersed.
\end{itemize}

Both points follow from simple calculations using lifts of
arcs to $\HH$.

\begin{lemma}\label{counting}
On $\xx$:
\begin{enumerate}
\item There are $2(p-1)$ arcs of $\lambda$-length $p$ incident to the cusp labelled  $\infty$;
\item These geodesics form $2(p-2)$ properly  immersed equilateral ideal triangles incident at $\infty$.
\end{enumerate}
\end{lemma}

\proof Suppose that $\alpha$ is a arc ending at the cusp labelled 
$\infty$.
It has a lift to $\HH$ which is a vertical line and after 
applying a tranformation in $\Gamma(2)$ of the 
form $z \mapsto z + 2n$ we may assume this lift has finite endpoint $\alpha^- \in [0,2]$.
By hypothesis $\alpha$ has  $\lambda$-length $p$ 
so that $\alpha^- = k/p$ for some $k$ co prime with $p$
and we have that
$$ k \in \{1,2,\ldots p-1, p+ 1, \ldots 2p -1 \}.$$
This proves the first point.

For the second point, begin by  noting that consecutive elements of this set of lifts
that is those ending at $k/p, (k+1)/p$ respectively, are sides of
an equilateral ideal triangle since
$$\det \begin{pmatrix}
k + 1 & k \\ p & p 
\end{pmatrix} = p.$$
It is easy to see that this is the only way to construct such
an  equilateral ideal triangle incident at $\infty$.
Thus there are $2(p-2)$ such triangles on $\xx$

\hfill $\Box$

\begin{lemma}\label{proper}
The projection of each of the triangles in the preceding lemma
is properly immersed.
\end{lemma}

\proof This follows from considerations of parity of numerator, denominator pairs.
By convention $\infty$ is $1/0$ and 
by hypothesis $p$ is odd so the denominators have different parities. So for $k$ co prime with $p$, 
$\infty$ and $k/p$ are in different $\Gamma(2)$ orbits.
Likewise the parity of $k$ and $k+1$ are always different
so that $k/p, (k+1)/p$ are in different $\Gamma(2)$ orbits.

\hfill $\Box$

\section{Concluding remarks}

We have proved Theorem \ref{eisenstein} by counting properly
immersed equilateral ideal triangles on the three punctured sphere.
Our proof is a simple application of the theory of hyperbolic
surfaces and the the fixed points of their automorphisms.
It is interesting to note that whilst our proof has similarities
with the proof of Theorem \ref{main} it is also quite different from
Heath-Brown's proof of Fermat's two squares theorem \cite{heath} as
it uses a different group of automorphisms and a different trick not
based on parity.

In further work in preparation we will investigate the relation  between values of
binary quadratic forms and $\lambda$-lengths of arcs
on other hyperbolic surfaces.

%\bibitem{Gui}
%L. Guillop\'e, Laurent(F-GREN-F)
%Sur la distribution des longueurs des g\'eod\'esiques ferm\'ees d'une surface compacte \`a bord totalement g\'eod\'esique. 
%Duke Math. J. 53 (1986), no. 3, 827-848. 

\end{document}